\documentclass[a4paper,12pt]{article}
\usepackage{amsmath}
\usepackage{amscd}
\usepackage[xdvi]{graphicx}
\usepackage[dvips]{color}
\usepackage{times}
\usepackage{txfonts}

\makeatletter

\@addtoreset{equation}{section}

\begin{document}

\begin{center}
\textbf{\LARGE{Continuous limit and the moments system for the globally coupled phase oscillators}}
\end{center}

\begin{center}
Institute of Mathematics for Industry, Kyushu University, Fukuoka,
819-0395, Japan

Hayato CHIBA \footnote{E mail address : chiba@imi.kyushu-u.ac.jp}
\end{center}
\begin{center}
Dec 1, 2011
\end{center}
\setlength{\baselineskip}{16pt}

\begin{center}
\textbf{Abstract}
\end{center}

The Kuramoto model, which describes synchronization phenomena,
is a system of ordinary differential equations on $N$-torus defined as coupled harmonic oscillators.
The order parameter is often used to measure the degree of synchronization.  
In this paper, the moments systems are introduced for both of the Kuramoto model and its continuous model.
It is shown that the moments systems for both systems take the same form.
This fact allows one to prove that the order parameter of the $N$-dimensional Kuramoto model converges to 
that of the continuous model as $N\to \infty$.


\section{Introduction}

Collective synchronization phenomena are observed in a variety of areas such as chemical reactions,
engineering circuits and biological populations~\cite{Pik}.
In order to investigate such a phenomenon, Kuramoto~\cite{Kura1, Kura2} proposed the system of ordinary differential equations
\begin{equation}
\frac{d\theta _i}{dt} 
= \omega _i + \frac{K}{N} \sum^N_{j=1} \sin (\theta _j - \theta _i),\,\, i= 1, \cdots  ,N,
\label{KMN}
\end{equation}
where $\theta _i \in [ 0, 2\pi )$ denotes the phase of an $i$-th oscillator on a circle,
$\omega _i\in \mathbf{R}$ denotes its natural frequency, $K>0$ is the coupling strength,
and where $N$ is the number of oscillators.
Eq.(\ref{KMN}) is derived by means of the averaging method from coupled dynamical systems having 
limit cycles, and now it is called the \textit{Kuramoto model}.

It is obvious that when $K=0$ or $K>0$ is sufficiently small,
$\theta _i(t)$ and $\theta _j(t)$ rotate on a circle at 
different velocities unless $\omega _i$ is equal to $\omega _j$.
On the other hand, if $K$ is sufficiently large, it is numerically observed that
some of oscillators or all of them tend to rotate at the same velocity on average, which is called the 
\textit{synchronization}~\cite{Pik,Str1}.
If $N$ is small, such a transition from de-synchronization to synchronization may be well revealed
by means of the bifurcation theory~\cite{ChiPa,Mai1,Mai2}.
However, if $N$ is large, it is difficult to investigate the transition from the view point of
the bifurcation theory and it is still far from understood.

\begin{figure}
\begin{center}
\includegraphics{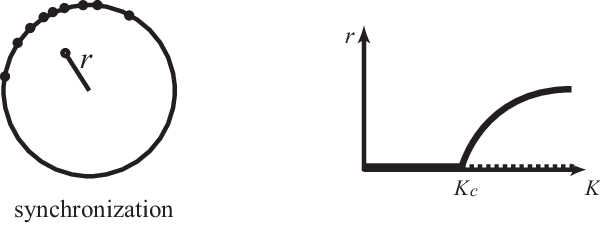}
\caption{A bifurcation diagram of the order parameter. \label{fig1}}
\end{center}
\end{figure}

In order to evaluate whether synchronization occurs or not, Kuramoto introduced
the \textit{order parameter} $r_N(t)$ by
\begin{equation}
r_N(t) := \frac{1}{N}\sum^N_{j=1} e^{\sqrt{-1} \theta _j(t)},
\label{order2}
\end{equation}
which gives the centroid of oscillators.
It seems that if synchronous state is formed, $|r_N(t)|$ takes a positive number, while
if de-synchronization is stable, $|r_N(t)|$ is zero on time average.
Indeed, based on some formal calculations, Kuramoto assumed a bifurcation diagram of the order parameter:
Suppose $N\to \infty$.
If $g(\omega )$, a distribution function for $\omega _i$'s,
is an even and unimodal function such that $g''(0)\neq 0$, then the bifurcation
diagram of $|r_\infty |$ is given as in Fig.\ref{fig1}. In other words, if the coupling strength $K$ is 
smaller than $K_c := 2/(\pi g(0))$, then $r_\infty \equiv 0$ is asymptotically stable.
If $K$ exceeds $K_c$, then a stable synchronous state emerges.
Near the transition point $K_c$, $|r_\infty |$ is of order $O((K- K_c)^{1/2})$.
See \cite{Str1} for Kuramoto's discussion.
In order to state his conjecture clearly, let us introduce the continuous model.

The infinite-dimensional version (the continuous model) of the Kuramoto model has been
well investigated to reveal a bifurcation diagram of the order parameter 
(see \cite{Ace,Bal,Chi,Cra,Mar,Mir2,Str1,Str2,Str3} and references therein).
The continuous model is defined as the equation of continuity of the form
\begin{equation}
\frac{\partial \rho_t}{\partial t} + \frac{\partial }{\partial \theta }(\rho_t v) = 0,
\end{equation}
where the unknown function $\rho_t = \rho_t (\theta , \omega )$ 
is a probability measure on $[0, 2\pi) \times \mathbf{R}$ parameterized by $t \in \mathbf{R}$.
Roughly speaking, $\rho_t (\theta , \omega )$ denotes a probability that
an oscillator having a natural frequency $\omega $ is placed at a position $\theta $.
See the next section for the definition of the vector field $v$.
The continuous version of the order parameter is defined to be 
\begin{equation}
r_\infty(t) := \int_{\mathbf{R}} \! \int^{2\pi}_{0} \! e^{\sqrt{-1}\theta } d\rho_t.
\end{equation}
Such a system is rather tractable because the order parameter for the infinite-dimensional version can be
constant in time, while the order parameter for the finite dimensional Kuramoto model is not constant
in general because solutions fluctuate due to effects of finiteness \cite{Str1}.
Recently, the Kuramoto's conjecture for the continuous model is rigorously proved by Chiba \cite{Chi};
The bifurcation diagram of the continuous version of the order parameter $r_\infty$ is given like as Fig.\ref{fig1}.

Now the questions arise : How close is the order parameter of the infinite-dimensional version
to that of the finite-dimensional Kuramoto model?
What is the influence of finite size effects?
This issue has been studied by many authors, see a reference paper \cite{Ace} by Acebron \textit{et al}.
In particular, Daido~\cite{Daido0} found the scaling law $|r_\infty - \langle r_N(t) \rangle| \sim (K_c - K)^{-1/2}N^{-1/2}$ for 
$K<K_c$ and $|r_\infty - \langle r_N(t) \rangle | \sim (K - K_c)^{-1/8}N^{-1/2}$ for $K>K_c$, although his analysis is not rigorous
from a view point of mathematics, where $r_\infty$ is assumed to be in a steady state, and
$\langle \,\, \rangle$ denotes the time average.  

In this paper, it is proved that the order parameter of the $N$-dimensional Kuramoto model converges to 
that of the continuous model in the sense of probability,
and their difference is of $O(N^{-1/2})$ as $N\to \infty$ for each $t$ (note that we do not take time average).
To prove this, the $(m,k)$-th moments are defined for both of the continuous model and the finite dimensional model.
In particular, $(0,1)$-th moment is the Kuramoto's order parameter. 
It is remarkable that both of the continuous model and the $N$-dimensional model
become the same evolution equation, called the \textit{moments system},
if they are rewritten by using the moments.
It means that any solutions of the continuous model and the $N$-dimensional model for any $N$
are embedded in the \textit{same} phase space of the moments system.
This fact allows us to measure the distance between solutions of the continuous model and that of 
the $N$-dimensional model in the same phase space.
These results and the central limit theorem prove that the difference between the order parameter of the 
$N$-dimensional model and that of the continuous model is of order $O(N^{-1/2})$ for each $t$,
provided that initial values and natural frequencies for the $N$-dimensional model are
independent and identically distributed according to a suitable probability measure.

The strategy of the proof is as follows:
Let $Z^m_k(t)$ and $\hat{Z}^m_k (t)$ be $(m,k)$-th moments of the continuous model and the $N$-dimensional model,
respectively (in particular, $Z^0_1 = r_\infty$ and $\hat{Z}^0_1 = r_N$).
Note that $Z^m_k(t)$ is determined by $\rho_t(\theta ,\omega )$ and $\hat{Z}^m_k(t)$ is determined by 
$\{ (\omega _j, \theta _j(t))\}_{j=1}^N$.
Let $h(\theta , \omega ) = \rho_0(\theta , \omega )$ be an initial measure for the continuous model.
Under a certain condition, there is a one-to-one correspondence between $\rho_t(\theta , \omega )$
and the set of moments $\{ Z^m_k(t)\}_{m,k}$ for each $t \geq 0$.
If initial values and natural frequencies for the $N$-dimensional model are
independent and identically distributed according to a measure $h(\theta , \omega )$,
the law of large number proves that $\hat{Z}^m_k (0) \to Z^m_k (0)$ as $N\to \infty$.
However, this argument is no longer applicable for a positive $t$ because 
$\theta _j(t)$'s are not independent and identically distributed when $t$ is positive.
Now we use the fact that \textit{$\{ Z^m_k(t)\}_{m,k}$ and $\{ \hat{Z}^m_k(t)\}_{m,k}$ are governed by the same 
differential equation called the moments system.}
Then, the continuity of solutions of the moments system with respect to initial conditions
immediately proves that if $\hat{Z}^m_k (t)$ and $Z^m_k (t)$ are sufficiently closed to one another for $t=0$,
the same is true for positive $t>0$. See the diagram below.

\begin{eqnarray*}
\begin{CD}
h(\theta , \omega ) @= \{ Z^m_k(0)\}_{m,k} @<N\to\infty<< \{ \hat{Z}^m_k(0)\}_{m,k} \\
@VVV @VVV @VVV \\
\rho_t(\theta , \omega ) @= \{ Z^m_k(t)\}_{m,k} @<N\to\infty<< \{ \hat{Z}^m_k(t)\}_{m,k}
\end{CD}
\end{eqnarray*}

Since the Kuramoto's conjecture for the continuous model is proved in \cite{Chi},
we obtain the following result as a corollary:
\begin{equation}
\lim_{t\to \infty}\lim_{N\to \infty} r_{N}(t) = \left\{ \begin{array}{ll}
0 & (0<K<K_c) \\
O(\sqrt{K-K_c}) & (K>K_c), \\
\end{array} \right.
\end{equation}
although behavior of another limit $\lim_{N\to \infty}\lim_{t\to \infty} r_{N}$ is still open.

More generally, a globally coupled phase oscillators defined to be
\begin{equation}
\frac{d\theta _i}{dt} 
= \omega _i + \frac{K}{N} \sum^N_{j=1} f(\theta _j - \theta _i),\,\, i= 1, \cdots  ,N,
\label{KDMN}
\end{equation}
is called the \textit{Kuramoto-Daido model} \cite{Daido4, Cra}, where the $2\pi$-periodic function 
$f : \mathbf{R} \to \mathbf{R}$ is called the \textit{coupling function}.
The results in this paper are easily extended to the Kuramoto-Daido model.


\section{Continuous model}

In this section, we introduce a continuous model of the Kuramoto model
and show existence, uniqueness and other properties of solutions of the model used in a later section.

Let us consider the Kuramoto model (\ref{KMN}). 
Following Kuramoto, we introduce the \textit{order parameter} $\hat{Z}^{0}_1$ by
\begin{equation}
\hat{Z}^{0}_1 (t) = \frac{1}{N} \sum^N_{j=1} e^{\sqrt{-1}\theta _j(t)}.
\label{order}
\end{equation}
The quantities $\hat{Z}^m_k$ will be defined in the next section.
By using it, Eq.(\ref{KMN}) is rewritten as
\begin{equation}
\frac{d \theta _i}{dt} = \omega _i 
      + \frac{K}{2\sqrt{-1}}(\hat{Z}^{0}_1 (t) e^{-\sqrt{-1}\theta _i}-\overline{\hat{Z}^{0}_1 (t)} e^{\sqrt{-1}\theta _i}),
\end{equation}
where $\overline{\hat{Z}^{0}_1}$ denotes the complex conjugate of $\hat{Z}^{0}_1$.
Motivated by these equations, we introduce a continuous model of the Kuramoto model,
which is an evolution equation of a probability measure $\rho_t = \rho_t (\theta , \omega )$ 
on $S^1 \times \mathbf{R}$ parameterized by $t \in \mathbf{R}$, as
\begin{eqnarray}
\left\{ \begin{array}{ll}
\displaystyle \frac{\partial \rho_t}{\partial t} + 
\frac{\partial }{\partial \theta }
\left( \Bigl(\omega  + \frac{K}{2\sqrt{-1}}(Z^{0}_1 (t) e^{-\sqrt{-1}\theta }-\overline{Z^{0}_1 (t)}
       e^{\sqrt{-1}\theta })\Bigr) \rho_t \right) = 0,  \\[0.3cm]
\displaystyle Z^0_1(t) := \int_{\mathbf{R}}
 \! \int^{2\pi}_{0} \! e^{\sqrt{-1}\theta } d\rho_t,  \\[0.3cm]
\rho_0 (\theta , \omega ) = h(\theta , \omega ),
\end{array} \right.
\label{conti}
\end{eqnarray}
where $h(\theta , \omega )$ is an initial measure.
The $Z^0_1(t)$ is a continuous version of $\hat{Z}^{0}_1 (t)$, and we also call it the \textit{order parameter}.
If we regard 
\begin{eqnarray*}
v_t := \omega  + \frac{K}{2\sqrt{-1}}(Z^{0}_1 (t) e^{-\sqrt{-1}\theta }-\overline{Z^{0}_1 (t)} e^{\sqrt{-1}\theta })
\end{eqnarray*}
as a velocity field, Eq.(\ref{conti}) provides an equation of continuity 
$\partial \rho_t / \partial t + \partial (\rho_t v_t)/ \partial \theta  = 0$ known in fluid dynamics.
It is easy to prove the law of conservation of mass:
\begin{equation}
\int_{\mathbf{R}} \! \int^{2\pi}_{0} \! \chi_{E} (\omega )  d\rho_t
 =  \int_{\mathbf{R}} \! \int^{2\pi}_{0} \! \chi_{E} (\omega )  dh =: g(E),
\label{g}
\end{equation}
where $E$ is any Borel set on $\mathbf{R}$ and $\chi_E (\omega )$ is the characteristic
function on $E$. A function $g$ defined as above gives a probability measure for natural frequencies
$\omega \in \mathbf{R}$ such that $\int_{\mathbf{R}} \! dg = 1$.
In particular $\int_{\mathbf{R}} \! \int^{2\pi}_{0} \! d\rho_t = 1$ if 
$\int_{\mathbf{R}} \! \int^{2\pi}_{0} \! dh = 1$.

By using the characteristic curve method, Eq.(\ref{conti}) is formally integrated as follows:
Consider the equation
\begin{eqnarray}
\frac{dx}{dt} &=& \omega + \frac{K}{2\sqrt{-1}}(Z^{0}_1 (t) e^{-\sqrt{-1}x}-\overline{Z^{0}_1 (t)} e^{\sqrt{-1}x}),\,\, x\in S^1,
\label{cha0} 
\end{eqnarray}
which defines a characteristic curve.
Let $x = x(t, s; \theta,\omega  )$ be a solution of Eq.(\ref{cha0}) 
satisfying the initial condition $x(s,s; \theta,\omega  ) = \theta $.
Along the characteristic curve $x(t):= x(t,0, \tilde{\theta }, \omega )$, 
$\rho_t (x(t),\omega )$ is differentiated as
\begin{eqnarray*}
\frac{d}{dt}\rho_t(x(t),\omega ) 
&=& \frac{\partial \rho_t}{\partial t}(x(t),\omega ) + \frac{\partial \rho_t}{\partial \theta }
\Bigl|_{\theta =x(t)} (\theta ,\omega ) \cdot \frac{dx}{dt}(t) \\
&=& \frac{\partial \rho_t}{\partial t}(x(t),\omega ) + \frac{\partial \rho_t}{\partial \theta }
\Bigl|_{\theta =x(t)} (\theta ,\omega ) \cdot v_t\bigl|_{\theta = x(t)} \\
&=& \frac{\partial \rho_t}{\partial t}(x(t),\omega ) + \frac{\partial}{\partial \theta }
\Bigl|_{\theta =x(t)} (v_t \rho_t )- \rho_t(x(t),\omega )\frac{\partial }{\partial \theta }
\Bigl|_{\theta =x(t)} v_t.
\end{eqnarray*}
Eq.(\ref{conti}) is used to yield
\begin{eqnarray*}
\frac{d}{dt}\rho_t(x(t),\omega ) 
&=& - \rho_t(x(t),\omega )\frac{\partial }{\partial \theta }\Bigl|_{\theta =x(t)} v_t \\
&=& \rho_t(x(t),\omega ) \cdot \frac{K}{2} 
\left( Z^{0}_1 (t) e^{-\sqrt{-1}x(t)}+\overline{Z^{0}_1 (t)} e^{\sqrt{-1}x(t)}\right).
\end{eqnarray*}
Hence, we obtain
\begin{eqnarray*}
\rho_t(x(t),\omega ) = h(\tilde{\theta }, \omega ) \cdot \exp \Bigl[ \frac{K}{2} \int^{t}_{0}\!
(Z^{0}_1 (s) e^{-\sqrt{-1}x(s)}+\overline{Z^{0}_1 (s)} e^{\sqrt{-1}x(s)} )ds \Bigr],
\end{eqnarray*}
which is true for any characteristic curve $x(t)= x(t,0, \tilde{\theta }, \omega )$.
Now we substitute $\tilde{\theta } = x(0,t; \theta ,\omega )$.
Due to the flow property, we have
\begin{eqnarray*}
x(s,0; x(0,t; \theta ,\omega ),\omega ) = x(s,t; \theta ,\omega ), \quad
x(t,0; x(0,t; \theta ,\omega ),\omega )  = \theta .
\end{eqnarray*}
Therefore, we obtain
\begin{equation}
\rho_t (\theta , \omega ) = h(x(0,t; \theta,\omega ), \omega ) \exp \Bigl[ 
\frac{K}{2} \int^t_{0} \! (Z^0_1 (s) e^{-\sqrt{-1} x(s,t; \theta,\omega ) }
                              + \overline{Z^0_1 (s)} e^{\sqrt{-1} x(s,t; \theta,\omega ) }) ds \Bigr],
\label{cha2}
\end{equation}
which gives a weak solution of (\ref{conti}).
By using Eq.(\ref{cha2}), it is easy to show the equality
\begin{equation}
\int_{\mathbf{R}} \! \int^{2\pi}_{0} \! a(\theta , \omega ) d\rho_t
 =  \int_{\mathbf{R}} \! \int^{2\pi}_{0} \! a(x(t,0; \theta,\omega ) , \omega ) dh,
\label{cha3}
\end{equation}
for any measurable function $a(\theta , \omega )$.
In particular, the order parameter $Z^0_1 (t)$ are rewritten as
\begin{equation}
Z^0_1 (t)
 =  \int_{\mathbf{R}} \! \int^{2\pi}_{0} \! e^{\sqrt{-1}x(t,0; \theta, \omega)} dh.
\label{cha4}
\end{equation}
Substituting it into Eqs.(\ref{cha0}), (\ref{cha2}), we obtain
\begin{equation}
\frac{d}{dt}x(t,s; \theta , \omega )
 = \omega  + K \int_{\mathbf{R}} \! \int^{2\pi}_{0} \! 
       \sin \Bigl( x(t,0;\theta ', \omega ') - x(t,s; \theta , \omega )\Bigr) dh(\theta ', \omega '),
\label{sol1}
\end{equation}
and
\begin{equation}
\rho_t (\theta , \omega ) = h(x(0,t; \theta,\omega ), \omega ) \exp \Bigl[ 
K \int^t_{0} \! ds \cdot \int_{\mathbf{R}} \! \int^{2\pi}_{0} \!\cos
\Bigl( x(s,0;\theta ', \omega ') - x(s,t; \theta , \omega ) \Bigr) dh(\theta ', \omega ') \Bigr]. 
\label{sol2}
\end{equation}
Even if $h(\theta , \omega )$ is not differentiable, we consider Eq.(\ref{sol2})
to be a weak solution of Eq.(\ref{conti}).
Indeed, even if $h$ and $\rho_t$ are not differentiable, the right hand side of (\ref{cha3}) is differentiable with respect to $t$
when $a(\theta , \omega )$ is differentiable.
\\[0.2cm]
\textbf{Theorem 2.1.} \, (i) There exists a unique weak solution $\rho_t$ of the initial value problem
(\ref{conti}) for any $t\geq 0$.
\\
(ii)  Solutions of (\ref{conti}) depend continuously on initial measures with respect to
the weak topology in the sense that for any numbers $T, \varepsilon >0$ and for any continuous
function $a(\theta , \omega )$ on $S^1\times \mathbf{R}$, there exist numbers $M(b)>0$ and
$\delta = \delta (T, \varepsilon , a) >0$ such that if initial measures $h_1, h_2$ satisfy
\begin{equation}
\left| \int_{\mathbf{R}} \! \int^{2\pi}_{0} \! b(\theta , \omega ) (dh_1 - dh_2) \right| < M(b)\delta ,
\label{cha5}
\end{equation}
for any continuous function $b(\theta , \omega )$,
then solutions $\rho_{t,1}$ and $\rho_{t,2}$
with $\rho_{0,1} = h_1$ and $\rho_{0,2} = h_2$ satisfy
\begin{equation}
\left| \int_{\mathbf{R}} \! \int^{2\pi}_{0} \! a(\theta , \omega ) (d\rho_{t,1} - d\rho_{t,2}) \right| < \varepsilon, 
\label{cha6}
\end{equation}
for $0 \leq t \leq T$. In particular if $a$ is Lipschitz continuous, then $\varepsilon  \sim O( \delta)$
as $ \delta \to 0$.
\\[0.2cm]
\textbf{Proof of (i).}\,  
It is sufficient to prove that the integro-ODE (\ref{sol1}) has a unique solution $x(t, s; \theta , \omega)$ satisfying
$x(s,s; \theta , \omega) = \theta $ for any $t,s \geq 0$ and $\theta \in S^1$.
Let us define a sequence $\{ x_n (t, 0; \theta, \omega )\}_{n=0}^{\infty}$ to be 
\begin{equation}
x_{n+1}(t,0;\theta, \omega ) = x_0(t,0;\theta, \omega )
 + K \int^t_{0} \! d\tau \cdot \int_{\mathbf{R}} \! \int^{2\pi}_{0} \!
f(x_{n}(\tau,0;\theta ', \omega') - x_{n}(\tau, 0;\theta , \omega)) dh(\theta ', \omega ')  
\end{equation}
and $x_0(t,0;\theta , \omega) = \theta  + \omega t$,
where $f(\theta ) = \sin \theta $ (since we prove the theorem for any $C^1$ function $f(\theta )$,
the theorem is also true for the continuous model for the Kuramoto-Daido model (\ref{KDMN})).
We estimate $|x_{n+1}(t,0;\theta , \omega)-x_{n}(t,0;\theta , \omega)|$ as
\begin{eqnarray*}
& & |x_{n+1}(t,0;\theta , \omega)-x_{n}(t,0;\theta , \omega)| \\
&\leq & K \! \int^t_{0} \! d\tau \cdot \!\int_{\mathbf{R}} \! \int^{2\pi}_{0} \!
\Bigl| f(x_{n}(\tau,0;\theta ', \omega') - x_{n}(\tau, 0;\theta , \omega)) \\
& & \quad \quad \quad     -  f(x_{n-1}(\tau,0;\theta ', \omega') - x_{n-1}(\tau, 0;\theta , \omega)) \Bigr| dh(\theta ', \omega' ) \\
&\leq & K L \!\int^t_{0} \! d\tau \cdot \!\int_{\mathbf{R}} \! \int^{2\pi}_{0} \!
\Bigl( |x_{n}(\tau,0;\theta ', \omega') - x_{n-1}(\tau, 0;\theta ', \omega')| \\
& & \quad \quad\quad    + |x_{n}(\tau,0;\theta , \omega) - x_{n-1}(\tau, 0;\theta , \omega)| \Bigr) dh(\theta ', \omega' ),
\end{eqnarray*}
where $L>0$ is the Lipschitz constant of the function  $f(\theta )$.
When $n=0$, we obtain
\begin{eqnarray*}
|x_{1}(t,0;\theta , \omega)-x_{0}(t,0;\theta , \omega)| \! \! 
&\leq & \!\!  K \!\!  \int^t_{0} \!\!  d\tau \cdot \!\! \int_{\mathbf{R}} \! \int^{2\pi}_{0} \!\! 
\left| f(x_{0}(\tau,0;\theta ', \omega') - x_{0}(\tau, 0;\theta , \omega)) \right| dh(\theta ', \omega' ) \\
&\leq & KMt,
\end{eqnarray*}
where $M = \max |f(\theta )|$.
Thus we can show by induction that
\begin{equation}
|x_{n}(t,0;\theta , \omega)-x_{n-1}(t,0;\theta , \omega)|  \leq 2^{n-1}L^{n-1}K^nM \frac{t^n}{n!}.
\end{equation}
This proves that $x_{n}(t,0;\theta , \omega)$ converges to a solution of the equation
\begin{eqnarray*}
\frac{dx}{dt}(t,0; \theta , \omega) 
  = \omega + K \int_{\mathbf{R}} \! 
          \int^{2\pi}_{0} \! f(x(t,0;\theta' , \omega')  - x(t,0;\theta , \omega) )dh(\theta ', \omega' ),
\end{eqnarray*}
as $n\to \infty$ for small $t\geq 0$.
Existence of global solutions are easily obtained by a standard way because 
the phase space $S^1$ is compact; that is,
solutions are extended for any $t>0$.
Uniqueness of solutions is also proved in a standard way and the detail is omitted.
With this $x(t,0; \theta , \omega)$, we define a sequence $\{ x_n (t, s; \theta , \omega)\}_{n=0}^{\infty}$ to be 
\begin{equation}
x_{n+1}(t,s;\theta , \omega) = x_0(t,s;\theta , \omega)
 + K \int^t_{0} \! d\tau \cdot \int_{\mathbf{R}} \! \int^{2\pi}_{0} \!
f(x(\tau,0;\theta ', \omega') - x_{n}(\tau, s;\theta , \omega)) dh(\theta ', \omega ')  
\end{equation}
and $x_0(t,s;\theta , \omega) = \theta  + \omega (t-s)$.
Then, existence and uniqueness of global solutions $x(t,s;\theta , \omega)$ is proved in the same way as above.
For this $x(t,s;\theta , \omega)$, Eq.(\ref{sol2}) gives a (weak) solution of Eq.(\ref{conti}).
\\[0.2cm]
\textbf{Proof of (ii).} \, Suppose that initial measures $h_1 ,h_2$ satisfy Eq.(\ref{cha5}).
Let $\rho_{t, 1}$ and $\rho_{t, 2}$ be solutions of Eq.(\ref{conti}) satisfying
$\rho_{0, 1} = h_1$ and $\rho_{0, 2} = h_2$. 
Let $x_i = x_i(t,0; \theta, \omega ),\,\, (i= 1,2)$ be solutions of 
\begin{equation}
\frac{dx_i}{dt}
 = \omega + K \int_{\mathbf{R}} \! \int^{2\pi}_{0} \! 
            f(x_i(t,0;\theta' , \omega')  - x_i(t,0;\theta , \omega) )dh_i(\theta ', \omega '),\,\, x_i\in S^1,
\end{equation}
satisfying $x_i(0,0; \theta , \omega) = \theta $, respectively.
Then we obtain
\begin{eqnarray}
& & \frac{d}{dt} \left( x_1 (t,0;\theta, \omega ) - x_2(t,0; \theta, \omega )\right) \nonumber \\
&=& K \int_{\mathbf{R}} \! \int^{2\pi}_{0} \! 
          f(x_1 (t,0;\theta', \omega' ) - x_1(t,0; \theta, \omega )) dh_1(\theta ', \omega' ) \nonumber \\
& &  -  K \int_{\mathbf{R}} \! \int^{2\pi}_{0} \! 
   f(x_2 (t,0;\theta', \omega' ) - x_2(t,0; \theta , \omega)) dh_2(\theta ', \omega ') \nonumber \\
&= & 
  K \int_{\mathbf{R}} \! \int^{2\pi}_{0} \! 
    f(x_1 (t,0;\theta' , \omega') - x_1(t,0; \theta, \omega )) (dh_1(\theta ', \omega' ) -dh_2(\theta ', \omega' )) \nonumber \\
& & + K \! \int_{\mathbf{R}} \! \int^{2\pi}_{0} \!\! \Bigl( 
f(x_1 (t,0;\theta' , \omega') - x_1(t,0; \theta, \omega ))\nonumber  \\
& & \quad \quad \quad    - f(x_2 (t,0;\theta', \omega' ) - x_2(t,0; \theta , \omega)) \Bigr) dh_2(\theta ', \omega' ).
\label{cha7}
\end{eqnarray}
Integrating it yields
\begin{eqnarray}
& & |x_1 (t,0;\theta, \omega ) - x_2(t,0; \theta, \omega )| \nonumber \\
&\leq & \int^t_{0} \! \Bigl| K \int_{\mathbf{R}} \! \int^{2\pi}_{0} \! 
    f(x_1 (s,0;\theta' , \omega') - x_1(s,0; \theta, \omega )) (dh_1(\theta ', \omega' ) -dh_2(\theta ', \omega' )) \Bigr| ds \nonumber\\
& &+ \int^t_{0} \!  K \! \int_{\mathbf{R}} \! \int^{2\pi}_{0} \!\! \Bigl| 
f(x_1 (s,0;\theta' , \omega') - x_1(s,0; \theta, \omega ))\nonumber  \\
& & \quad \quad \quad    - f(x_2 (s,0;\theta', \omega' ) - x_2(s,0; \theta , \omega)) \Bigr| dh_2(\theta ', \omega' ) ds \nonumber\\
& \leq & KM\delta t + KL \int^t_{0} \!  \! \int_{\mathbf{R}} \! \int^{2\pi}_{0} \!\! \Bigl( 
|x_1 (s,0;\theta' , \omega') - x_2(s,0; \theta', \omega' )| \nonumber  \\
& & \quad \quad \quad    + |x_1 (s,0;\theta, \omega ) - x_2(s,0; \theta , \omega)| \Bigr) dh_2(\theta ', \omega' ) ds,
\label{cha7b}
\end{eqnarray}
where $L$ is the Lipschitz constant of $f$ and $M$ is a constant arising from Eq.(\ref{cha5}).
If we put 
\begin{eqnarray*}
F(t, \theta , \omega ) = |x_1 (t,0;\theta, \omega ) - x_2(t,0; \theta, \omega )|,
\end{eqnarray*}
then (\ref{cha7b}) provides
\begin{eqnarray*}
\int_{\mathbf{R}} \! \int^{2\pi}_{0} \! F(t, \theta , \omega ) dh_2 \leq KM\delta t
           + 2KL \int^t_{0} \! \int_{\mathbf{R}} \! \int^{2\pi}_{0} \! F(s,\theta ,\omega ) dh_2 ds . 
\end{eqnarray*}
Now the Gronwall inequality proves
\begin{eqnarray*}
\int_{\mathbf{R}} \! \int^{2\pi}_{0} \! F(t, \theta , \omega ) dh_2 \leq \frac{M\delta }{2L}(e^{2KLt} - 1).
\end{eqnarray*}
Substituting it into (\ref{cha7b}) yields
\begin{eqnarray*}
& & |x_1 (t,0;\theta, \omega ) - x_2(t,0; \theta, \omega )| \nonumber \\
& \leq & KM\delta t + \frac{KM\delta }{2}\int^t_{0} \! (e^{2KLs} - 1) ds
 + KL\int^t_{0} \! |x_1 (s,0;\theta, \omega ) - x_2(s, 0; \theta, \omega )|ds \nonumber \\
& \leq & \frac{KM\delta t}{2} + \frac{M\delta }{4L}(e^{2KLt} - 1)
+ KL \int^t_{0} \! |x_1 (s,0;\theta, \omega ) - x_2(s, 0; \theta, \omega )|ds \nonumber.
\end{eqnarray*}
The Gronwall's inequality is applied again to obtain
\begin{equation}
|x_1 (t,0;\theta, \omega ) - x_2(t,0; \theta, \omega )| \leq \frac{M\delta }{2L}(e^{2KLt} - 1).
\label{cha8}
\end{equation}
Finally, the left hand side of Eq.(\ref{cha6}) is estimated as
\begin{eqnarray}
& & \left| \int_{\mathbf{R}} \! \int^{2\pi}_{0} \! a(\theta , \omega ) (d\rho_{t,1} - d\rho_{t,2}) \right| \nonumber \\
&= & \left| \int_{\mathbf{R}} \! \int^{2\pi}_{0} \!  \Bigl( a(x_1(t,0; \theta, \omega ) , \omega ) dh_1
           - a(x_2(t,0;\theta, \omega ), \omega )dh_2 \Bigr) \right| \nonumber \\
& \leq & \int_{\mathbf{R}} \! \int^{2\pi}_{0} \!\! 
         \left| a(x_1(t,0; \theta, \omega ) , \omega ) - a(x_2(t,0; \theta , \omega) , \omega ) \right| dh_2 \nonumber \\
& & \quad \quad    + \left| \int_{\mathbf{R}} \!\! \int^{2\pi}_{0} \!\!\!  a(x_1(t,0;\theta , \omega),\omega )(dh_1 - dh_2) \right| .
\label{cha9}
\end{eqnarray}
Since $a(\theta ,\omega ) $ is continuous and since Eq.(\ref{cha8}) holds, 
the first term in the right hand side of the above is less than $\varepsilon /2$ for $0 \leq t\leq T$ if $\delta$ 
is sufficiently small. The second term is also less than $\varepsilon /2$ if $\delta$
is sufficiently small because of Eq.(\ref{cha5}).
This proves Eq.(\ref{cha6}).
It is easy to see by Eq.(\ref{cha9}) that if $a(\theta ,\omega ) $ is Lipschitz continuous, then $\varepsilon $ 
is of order $O(\delta )$.  \hfill $\blacksquare$


\section{Moments system}

In this section, we introduce a moments system to transform the finite-dimensional Kuramoto model
(\ref{KMN}) and its continuous model (\ref{conti}) into the same system.
We prove by using the moments system that the order parameter (\ref{order2}) for the Kuramoto model
converges to the order parameter $Z^0_1(t)$ for the continuous model as $N\to \infty$ under 
appropriate assumptions.
\\

For a given probability measure $h(\theta , \omega )$ on $S^1 \times \mathbf{R}$,
suppose that absolute moments
\begin{equation}
M^n_k := \int_{\mathbf{R}} \! \int^{2\pi}_{0} \! |\omega ^n e^{\sqrt{-1}k \theta }| dh 
\label{m1} 
\end{equation}
exist for $k=0,\pm 1, \cdots $ and $m=0,1,\cdots $. Then, the moments $m^n_k$ are defined to be
\begin{equation}
m^n_k := \int_{\mathbf{R}} \! \int^{2\pi}_{0} \! \omega ^n e^{\sqrt{-1}k \theta } dh.
\label{m2} 
\end{equation}
Conversely, if there exists a unique probability measure  $h$ for a given sequence of numbers $\{m^n_k\}_{n,k}$
such that Eq.(\ref{m2}) holds, then $h$ is called \textit{M-determinate}.
In this case, we also say that moments $\{m^n_k\}_{n,k}$ is M-determinate.
Many conditions for which $h$ is M-determinate have been well studied
as the moment problem \cite{Akh,Sho,Fro,Sim}.
For example, one of the most convenient conditions is that if $h$ has all absolute moments
$M^n_k$ and they satisfy $\sum^\infty_{n=1} (M^n_0 + 1)^{-1/n} = \infty$ (\textit{Carleman's condition}),
then $h$ is M-determinate.
\\[0.2cm]
\textbf{Example 3.1.}\, If $h$ has compact support, then $h$ is M-determinate.
Suppose that $h$ has a probability density function of the form $\hat{h}(\theta ) \hat{g}(\omega )$.
If $\hat{g}(\omega )$ is the Gaussian distribution, then $h$ is M-determinate.
If $\hat{g}(\omega ) = 1/(\pi (1 + \omega ^2))$ is the Lorentzian distribution, $h$ is 
\textit{not} M-determinate because $h$ does not have all moments.
\\

In what follows, we suppose that an initial measure $h(\theta ,\omega )$ for the initial value problem 
(\ref{conti}) has all absolute moments and is M-determinate.
Recall that a probability measure $g$ for the natural frequency $\omega $
is defined through Eq.(\ref{g}). Since $h$ has absolute moments
\begin{equation}
M^n_k = \int_{\mathbf{R}} \! \int^{2\pi}_{0} \! |\omega ^n e^{\sqrt{-1}k\theta }| dh = \int_{\mathbf{R}} \! |\omega |^n dg<\infty, 
\label{m3} 
\end{equation}
$g$ also has all moments $\mu_n := \int_{\mathbf{R}} \! \omega ^n dg,\, n = 0,1,2,\cdots $.
Consider the Lebesgue space $L^2 (\mathbf{R}, dg)$.
Since all moments $\mu_m$ of $g$ exist, we can construct a complete orthonormal system
$\{P_m(\omega )\}_{m=0}^\infty$ on $L^2 (\mathbf{R}, dg)$,
by using the Gram-Schmidt orthogonalization from $\{\omega ^m\}_{m=0}^\infty$, such that
\begin{equation}
( P_n, P_m ) = \int_{\mathbf{R}} \! P_n(\omega ) P_m(\omega ) dg = \left\{ \begin{array}{ll}
1 & (n=m) \\
0 &  (n\neq m), \\
\end{array} \right.
\label{m4} 
\end{equation}
where $(\,\, , \,\, )$ denotes the inner product on $L^2 (\mathbf{R}, dg)$ and $P_n(\omega )$ is a polynomial
of degree $n$. In particular, $P_0(\omega ) \equiv 1$.
It is well known that $P_n(\omega )$ satisfies the relation
\begin{equation}
\omega P_n(\omega ) = b_n P_{n+1} (\omega ) + a_n P_n(\omega ) + b_{n-1} P_{n-1}(\omega )
\label{m5} 
\end{equation}
for $n = 0,1,2 \cdots $, where $a_n$ and $b_n$ are real constants determined by $g$.
The matrix $\mathcal{M}$ defined as
\begin{equation}
\mathcal{M} = \left(
\begin{array}{@{\,}ccccc@{\,}}
a_0 & b_0 & & & \\
b_0 & a_1 & b_1 & & \\
 & b_1& a_2 & b_2 & \\
 & & & \ddots & \\
 & & & & \\
\end{array}
\right)
\label{m6} 
\end{equation}
is called the \textit{Jacobi matrix} for $g$.
Eq.(\ref{m5}) shows that the Jacobi matrix gives the $l^2 (\mathbf{Z}_{\geq 0})$ representation
of the multiplication operator
\begin{equation}
\mathcal{M} : p(\omega ) \mapsto \omega p(\omega )
\label{m7} 
\end{equation}
on $L^2 (\mathbf{R}, dg)$, where 
$l^2 (\mathbf{Z}_{\geq 0}) = \{ \{x_n\}_{n=0}^\infty \, | \, \sum^\infty_{n=0} |x_n|^2 < \infty\}$.
\\

If an initial measure $h$ is M-determinate,
so is a solution $\rho_t$ of the continuous model (\ref{conti}) because of Eq.(\ref{sol2}).
Let us define the \textit{$(m,k)$-th moments} $Z^m_k$ for $\rho_t$ to be
\begin{equation}
Z^m_k (t) = \int_{\mathbf{R}} \! \int^{2\pi}_{0} \! P_m(\omega ) e^{\sqrt{-1} k \theta } d\rho_t,
\label{m8} 
\end{equation}
for $m = 0,1, 2 \cdots $ and $k= 0, \pm 1, \pm 2, \cdots $.
In particular $Z^0_1(t)$ is the order parameter given in Eq.(\ref{conti}),
and $Z^m_{-k}(t) = \overline{Z^m_k(t)}$.
Note that
\begin{equation}
Z^m_0(t) = 
\left\{ \begin{array}{ll}
1 & (m=0) \\
0 &  (m\neq 0)\\
\end{array} \right.
\label{m9} 
\end{equation}
are constants. It is easy to verify that
\begin{equation}
|Z^m_k(t)| \leq 1
\label{m10} 
\end{equation}
by using the Schwarz inequality.
By using Eq.(\ref{cha3}), an evolution equation for $Z^m_k (t)$ is derived as
\begin{eqnarray}
\frac{dZ^m_k}{dt} \!\!\! &=& \frac{\partial }{\partial t}
          \int_{\mathbf{R}} \! \int^{2\pi}_{0} \! P_m(\omega ) e^{\sqrt{-1} k x(t,0;\theta, \omega )}dh \nonumber \\
&=&  \int_{\mathbf{R}} \! \int^{2\pi}_{0} \! P_m(\omega ) \sqrt{-1}k \frac{\partial x}{\partial t}
                 (t,0;\theta, \omega ) e^{\sqrt{-1} k x(t,0;\theta , \omega)} dh \nonumber \\
&=& \!\!\!  \int_{\mathbf{R}} \! \int^{2\pi}_{0} \!\! P_m(\omega ) \sqrt{-1}k 
     \left( \omega + \frac{K}{2\sqrt{-1}} ( Z^0_1 (t) e^{- \sqrt{-1} x(t,0;\theta , \omega)} - Z^0_{-1} (t) e^{ \sqrt{-1} x(t,0;\theta , \omega)})
              \right) e^{\sqrt{-1} k x(t,0;\theta, \omega )} dh \nonumber \\
&=&\sqrt{-1}k \int_{\mathbf{R}} \! \int^{2\pi}_{0} \! \omega  P_m(\omega ) e^{\sqrt{-1} k x(t,0;\theta, \omega )}dh \nonumber \\
& & \quad + \frac{kK}{2} \int_{\mathbf{R}} \! \int^{2\pi}_{0} \! P_m(\omega )
         (Z^0_1(t) e^{ \sqrt{-1} (k-1) x(t,0;\theta, \omega )} - Z^0_{-1}(t) e^{ \sqrt{-1} (k+1) x(t,0;\theta, \omega )}) dh \nonumber \\
&=& \sqrt{-1}k \int_{\mathbf{R}} \! \int^{2\pi}_{0} \! \Bigl( b_m P_{m+1} (\omega ) 
           + a_m P_m(\omega ) + b_{m-1} P_{m-1}(\omega ) \Bigr) e^{\sqrt{-1} k \theta } d\rho_t \nonumber \\
& & \quad + \frac{kK}{2} \int_{\mathbf{R}} \! \int^{2\pi}_{0} \! P_m(\omega )
      (Z^0_1(t) e^{ \sqrt{-1} (k-1) \theta} - Z^0_{-1}(t) e^{ \sqrt{-1} (k+1) \theta }) d\rho_t \nonumber \\
&=& \sqrt{-1}k \left( b_m Z^{m+1}_k + a_m Z^m_k + b_{m-1} Z^{m-1}_{k} \right) 
          + \frac{kK}{2} (Z^0_1 Z^m_{k-1} -Z^0_{-1} Z^m_{k+1}).
\label{m11} 
\end{eqnarray}
Put $Z_k = (Z^0_k, Z^1_k, Z^2_k, \cdots )^T$, where $T$ denotes the transpose.
Define the Jacobi matrix $\mathcal{M}$ and the projection matrix $\mathcal{P}$ to be Eq.(\ref{m6}) 
and 
\begin{equation}
\mathcal{P} = \left(
\begin{array}{@{\,}cccc@{\,}}
1& 0 & \cdots & \\
0 &0 & & \\
\vdots& & \ddots  & \\
& & & 
\end{array}
\right),
\label{m12} 
\end{equation}
respectively. Then, Eq.(\ref{m11}) is rewritten as
\begin{eqnarray}
\frac{d}{dt} \left(
\begin{array}{@{\,}c@{\,}}
Z_1 \\
Z_2 \\
Z_3 \\
 \vdots
\end{array}
\right) =  \left(
\begin{array}{@{\,}cccc@{\,}}
\! \displaystyle \sqrt{-1} \mathcal{M} + \frac{K}{2} \mathcal{P} \!\! & & \\
& \!\!\!\! 2\sqrt{-1} \mathcal{M}  \!\!\!\! & & \\
& & \!\!\!\! 3\sqrt{-1} \mathcal{M} \!\!\!\!& \\
& & & \ddots
\end{array}
\right) \left(
\begin{array}{@{\,}c@{\,}}
Z_1 \\
Z_2 \\
Z_3 \\
\vdots
\end{array}
\right) + \frac{K}{2} \left(
\begin{array}{@{\,}c@{\,}}
-Z^0_{-1}\, Z_{2} \\
2(Z^0_{1}\, Z_{1} - Z^0_{-1}Z_3) \\
3(Z^0_{1}\, Z_{2} - Z^0_{-1}Z_4) \\
\vdots
\end{array}
\right) .
\label{m13}
\end{eqnarray}
Note that equations for $Z_{-1}, Z_{-2}, \cdots $ are omitted because $Z_{-k} = \overline{Z}_k$.
The first term is a linear term and the second is a nonlinear term.
We call Eq.(\ref{m11}) or Eq.(\ref{m13}) the \textit{moments system}.
The dynamics of the system is investigated in \cite{Chi}.

Let $M_D$ be the set of M-determinate sequences $\{Z^m_k\}_{m,k}$ in the sense that if 
$\{Z^m_k\}_{m,k}\in M_D$, then there exists a unique measure $h$ on $S^1 \times \mathbf{R}$
such that $Z^m_k = \int_{\mathbf{R}} \! \int^{2\pi}_{0} \! P_m(\omega ) e^{\sqrt{-1}k \theta } dh $.
Since there is a one-to-one correspondence between elements of $M_D$ and M-determinate measures,
Thm.2.1 is restated as follows.
\\[0.2cm]
\textbf{Theorem 3.2.}\, (i) There exists a unique
solution $\{Z^m_k (t) \}_{m,k} \in M_D$ of the moments system if 
an initial condition $\{Z^m_k (0) \}_{m,k} $ is in $M_D$.
\\
(ii) Let $\{ Z^m_k(t)\}$ and $\{ \widetilde{Z}^m_k (t)\}$ be solutions of the moments system
with initial conditions $\{ Z^m_k(0)\}$, $\{ \widetilde{Z}^m_k (0)\} \in M_D$, respectively.
For any positive numbers $T$ and $\varepsilon $, there exist positive numbers $C_{m,k}$ and 
$\delta = \delta (T, \varepsilon )$ such that if 
\begin{equation}
|Z^m_k (0) - \widetilde{Z}^m_k (0)| < C_{m,k} \delta,
\label{m14}
\end{equation}
for any $m,k$, then the inequality
\begin{equation}
|Z^m_k (t) - \widetilde{Z}^m_k (t)| < \varepsilon 
\label{m15}
\end{equation}
holds for $0 \leq t \leq T$.
In particular $\varepsilon \sim O(\delta )$ as $\delta \to 0$.
\\

For the $N$-dimensional Kuramoto model (\ref{KMN}), 
we define the \textit{$(m,k)$-th moments} to be
\begin{equation}
\hat{Z}^m_k (t) := \frac{1}{N} \sum^N_{j=1} P_m(\omega _j) e^{\sqrt{-1} k \theta _j(t)},
\label{m16}
\end{equation}
for $m=0,1,2,\cdots $ and $k = 0, \pm 1, \pm 2, \cdots $.
In particular $\hat{Z}^0_1$ is the order parameter defined in Eq.(\ref{order}).
By using Eq.(\ref{KMN}), it is easy to verify that $\hat{Z}^m_k (t)$'s 
satisfy a system of differential equations
\begin{equation}
\frac{d\hat{Z}^m_k (t)}{dt}
 = \sqrt{-1}k \left( b_m \hat{Z}^{m+1}_k + a_m \hat{Z}^m_k + b_{m-1} \hat{Z}^{m-1}_{k} \right) 
          + \frac{kK}{2} (\hat{Z}^0_1 \hat{Z}^m_{k-1} - \hat{Z}^0_{-1} \hat{Z}^m_{k+1}).
\label{m17}
\end{equation}
It is remarkable that Eq.(\ref{m17}) has the same form as Eq.(\ref{m11}).
This means that all solutions of the Kuramoto model for any $N$
are embedded in the phase space of the moments system (\ref{m11}).
This fact allows us to prove Theorem 3.3 below.
Originally the moments system for the Kuramoto model was introduced by Perez and Ritort~\cite{Per},
although their definition of the moments is $H^m_k := 1/N \cdot \sum^N_{j=1} \omega ^m_j e^{\sqrt{-1}k \theta _j(t)}$.
Since we adopt orthogonal polynomials $\{ P_m(\omega )\}_{m=0}^{\infty}$ to define moments (\ref{m16}),
our moments system is more suitable for mathematical analysis.
\\

Now we are in a position to show the main theorem in this paper,
which states that differences between moments $Z^m_k (t)$ and  $\hat{Z}^m_k (t)$
are of $O(1/\sqrt{N})$ and thus the continuous model (\ref{conti}) is proper to
investigate the Kuramoto model (\ref{KMN}) for large $N$.
\\[0.2cm]
\textbf{Theorem 3.3.} \, Let $\rho_t$ be a solution of the continuous model (\ref{conti}) 
such that an initial measure $h(\theta , \omega )$ is M-determinate.
Suppose that for the $N$-dimensional Kuramoto model (\ref{KMN}),
pairs $(\theta _j(0), \omega _j)$ of 
initial values $\theta _j(0),\, j= 1, \cdots  ,N$ and natural frequencies $\omega _j,\, j= 1, \cdots  ,N$
are independent and identically distributed according to the probability measure $h(\theta , \omega )$.
Then, moments $Z^m_k (t)$ and $\hat{Z}^m_k (t)$ defined by Eqs.(\ref{m8}) and (\ref{m16})
satisfy
\begin{equation}
|Z^m_k (t) - \hat{Z}^m_k (t)| \to 0 ,\,\,\,\, a.s. \,\, (N \to \infty),
\label{m18}
\end{equation}
for any $m,k$ and $t$ ($a.s.$ denotes `` almost surely").
Further, for any positive number $\delta $, there exists a number $C = C(m,k,t, \delta) > 0$
such that
\begin{equation}
P (\, |Z^m_k (t) - \hat{Z}^m_k (t)|< C/ \sqrt{N} \,) \to 1-\delta,\,\, (N\to \infty),
\label{m19}
\end{equation}
where $P(A)$ is the probability that an event $A$ will occur.
\\[0.2cm]
\textbf{Proof.} 
Since $\omega _j$'s and $\theta _j(0)$'s are independent and identically distributed,
the average of $\hat{Z}^m_k (0)$ is calculated as
\begin{eqnarray}
E[\hat{Z}^m_k (0)] &=& E[ \frac{1}{N} \sum^N_{j=1} P_m(\omega _j) e^{\sqrt{-1} k \theta _j(0)}] \nonumber \\
&=& E[P_m(\omega _j) e^{\sqrt{-1} k \theta _j(0)}] \nonumber \\
&=& \int_{\mathbf{R}} \! \int^{2\pi}_{0} \! P_m(\omega ) e^{\sqrt{-1} k \theta } dh = Z^m_k(0). 
\label{m20}
\end{eqnarray}
Thus Eqs.(\ref{m18}) and (\ref{m19}) for $t=0$ immediately follow from the strong law of large number
and the central limit theorem, respectively.
Note that the strong law of large number and the central limit theorem are no longer applicable for $t >0$
because $\theta _j(t)$'s are not independent and identically distributed when $t$ is positive.
However, since $Z^m_k (t)$ and $\hat{Z}^m_k (t)$ satisfy the same moments system,
and since solutions of the moments system are continuous with respect to initial values (Thm.3.2 (ii)),
Eqs.(\ref{m18}),(\ref{m19}) hold for each positive $t$ if they are true for $t=0$;
if initial states satisfy $|Z^m_k (0) - \hat{Z}^m_k (0)| \to 0$ as $N\to \infty$,
then $|Z^m_k (t) - \hat{Z}^m_k (t)| \to 0$ holds for any $t > 0$, and if they satisfy
$|Z^m_k (0) - \hat{Z}^m_k (0)|< C_0/ \sqrt{N}$ for a positive constant $C_0$,
then $|Z^m_k (t) - \hat{Z}^m_k (t)|< C_t/ \sqrt{N}$ holds for some $C_t >0$.
 \hfill $\blacksquare$
\\

\textbf{Acknowledgements}

This work was supported by Grant-in-Aid for Young Scientists (B), No.22740069 from MEXT Japan.


\end{document}